\documentclass[
final
]{dmtcs-episciences}

\usepackage[utf8]{inputenc}
\usepackage{subfigure}
\usepackage{tikz}
\usepackage{amssymb,amscd,amsmath,amsfonts,amsthm}
\usepackage{enumerate}
\usepackage{multicol}
\usepackage{hyperref}

\newtheorem{theorem}{Theorem}

\newtheorem{corollary}[theorem]{Corollary}

\newtheorem{lemma}[theorem]{Lemma}

\newtheorem{observation}[theorem]{Observation}
\newtheorem{problem}{Problem}
\newtheorem{proposition}[theorem]{Proposition}
\newtheorem{remark}[theorem]{Remark}

\newcommand{\ra}{{\rm a}}
\newcommand{\rb}{{\rm b}}
\newcommand{\cF}{{\cal F}}
\newcommand{\cS}{{\cal S}}

\usepackage[round]{natbib}

\author[Bo\v stjan Bre\v sar et al.]{Bo\v stjan Bre\v sar\affiliationmark{1,2}\thanks{Supported in part by the Slovenian Research and Innovation Agency (ARIS) under the grants P1-0297, N1-0285, N1-0431 and J1-4008.}
  \and Csilla Bujt\'{a}s\affiliationmark{2,3}\thanks{Supported in part by the Slovenian Research and Innovation Agency (ARIS) under the grants P1-0297 and N1-0355.}
  \and Pakanun Dokyeesun\affiliationmark{2}\thanks{Supported in part by the Slovenian Research and Innovation Agency (ARIS) under the grant P1-0297.}
  \and Tanja Dravec\affiliationmark{1,2}\thanks{Supported in part by the Slovenian Research and Innovation Agency (ARIS) under the grants P1-0297 and N1-0431.}}
\title[Thresholds for the biased Maker-Breaker domination games]{Thresholds for the biased Maker-Breaker domination games}
\affiliation{
  Faculty of Natural Sciences and Mathematics, University of Maribor, Maribor, Slovenia\\
  Institute of Mathematics, Physics and Mechanics, Ljubljana, Slovenia\\
  Faculty of Mathematics and Physics, University of Ljubljana, Ljubljana, Slovenia}
\begin{document}
% This is only used if you are compiling for a volume before vol 25
% \publicationdetails{VOL}{2015}{ISS}{NUM}{SUBM}
% This is the new form of collecting the data, starting with vol 25
\publicationdata
{vol. 27:3}
{2025}
{19}
{10.46298/dmtcs.15392}
%{1998-10-14; 1998-10-14; 2002-07-19; 2014-02-05; 2015-09-09; 2022-12-25}
%{2022-12-3}
{2025-03-18; 2025-03-18; 2025-09-04}
{2025-10-14}
\maketitle
\begin{abstract}
  In the $(a,b)$-biased Maker-Breaker domination game, two players alternately select unplayed vertices in a graph $G$ such that Dominator selects $a$ and Staller selects $b$ vertices per move. Dominator wins if the vertices he selected during the game form a dominating set of $G$, while Staller wins if she can prevent Dominator from achieving this goal. Given a positive integer $b$, Dominator's threshold, ${\rm a}_b$, is the minimum $a$ such that Dominator wins the $(a,b)$-biased game on $G$ when he starts the game. Similarly, ${\rm a}'_b$ denotes the minimum $a$ such that Dominator wins when Staller starts the $(a,b)$-biased game. Staller's thresholds, ${\rm b}_a$ and ${\rm b}'_a$, are defined analogously. It is proved that Staller wins the $(k-1,k)$-biased games in a graph $G$ if its order is sufficiently large with respect to a function of $k$ and the maximum degree of $G$. Along the way, the $\ell$-local domination number of a graph is introduced. This new parameter is proved to bound Dominator's thresholds ${\rm a}_\ell$ and ${\rm a}_\ell'$ from above. As a consequence, ${\rm a}_1'(G)\le 2$ holds for every claw-free graph $G$. More specific results are obtained for thresholds in line graphs and Cartesian grids. 
Based on the concept of $[1,k]$-factor of a graph $G$, we introduce the star partition width $\sigma(G)$ of $G$, and prove that ${\rm a}_1'(G)\le \sigma(G)$ holds  for any nontrivial graph $G$, while ${\rm a}_1'(G)=\sigma(G)$ if $G$ is a tree.
\end{abstract}

\section{Introduction}
\label{sec: intro}

The \emph{Maker-Breaker domination game} (\emph{MBD game} for short) is played on a graph $G$ by two players, called Dominator and Staller. The players alternately select an unplayed vertex in $G$. Dominator wins the game if he can claim a dominating set of $G$ while Staller wins if she can prevent Dominator from claiming a dominating set. In other words, Staller wins if she claims a closed neighborhood of a vertex in $G$.
There are two versions of the MBD game based on who starts it. 
A Dominator-start MBD game will also be called a \emph{D-game}, while \emph{S-game} is a short name for the  Staller-start MBD game. 

Maker-Breaker games are among the most important positional games; see~\cite{beck, beck-1981,  erdos-1973, hefetz-2014} for more detail.  In many cases, Maker-Breaker games have been extended by considering biased variations of the original games (see e.g., \cite{chvatal-1978,clem-2018}), in which the players, on their turn, select an arbitrary but fixed number of vertices.  Domination games (see~\cite{bresar-2010, book-2021}) and Maker-Breaker domination games (see~\cite{bdg-2025, bujtas-2021, bujtas-2023, dokyeesun-2024}) have been extensively investigated in recent years, and very recently, the study of their biased versions has been initiated in~\cite{bagdas-2024}.  Let $a,b$ be two positive integers and $G$ be a graph. The \emph{$(a,b)$-biased Maker-Breaker domination game} is the same as the Maker-Breaker game on $G$, except that in the biased game Dominator chooses $a$ unplayed vertices per move and Staller chooses $b$ unplayed vertices per move. The integers $a,b$ are called the \emph{bias} of Dominator and Staller, respectively.
 In the last move of the game, if there are fewer unplayed vertices than his/her bias, a player selects all remaining unplayed vertices.
The Maker-Breaker domination game is a special case of the  $(a,b)$-biased Maker-Breaker domination game where $a=b=1$.

Assuming that the $(a,b)$-biased MBD game is played on a graph $G$ and both players are playing optimally according to their goals, let $W(G,a,b)$ denote the winner in the D-game, while $W'(G,a,b)$ is the winner in the S-game. Abbreviating the names of Dominator and Staller to $D$ and $S$, we thus have $W(G,a,b)\in \{D,S\}$ and $W'(G,a,b)\in \{D,S\}$. 
The $(a,b)$-biased MBD game when Dominator is the first player is called the \emph{$(a,b)$-biased D-game}. Similarly, the \emph{$(a,b)$-biased S-game} is the $(a,b)$-biased MBD game when Staller starts. These games were first studied in a recent manuscript of~\cite{bagdas-2024}, which focuses on the $(a,1)$-biased games and also studies the number of moves needed to win in such games.   

Let $G$ be a graph. When $b$ is fixed, then we are interested in the minimum $a$ such that $W(G,a,b)=D$, respectively $W'(G,a,b)=D$, which yields the so-called Dominator's thresholds, denoted by $\ra_b(G)$ and $\ra_b'(G)$, respectively.  In a similar way, one defines Staller's thresholds, $\rb_a(G)$ and $\rb_a'(G)$.
If $G$ is a graph of infinite order but the vertex degrees are finite, the Maker-Breaker domination game on $G$ is defined analogously. In this case, Staller wins the game if she can claim an entire closed neighborhood of a vertex, and Dominator wins if he can prevent Staller from claiming a closed neighborhood. The other concepts as the $(a,b)$-biased D-game, S-game, and the thresholds $\ra_b(G)$, $\ra_b'(G)$, $\rb_a(G)$, $\rb_a'(G)$ can also be extended to infinite graphs with finite vertex degrees. In Section~\ref{sec:preliminaries}, we give formal definitions of these concepts. 

\paragraph{Structure of the paper.} The main goal of this paper is to obtain exact values or bounds for the thresholds in $(a,b)$-biased MBD games on graphs. In Section~\ref{sec:large-order}, we prove that for every graph $G$ of a sufficiently large order, which depends on a function of $k$ and maximum degree $\Delta(G)$, Staller wins the $(k-1,k)$-biased game. From this we infer that for graphs $G$ whose order is bounded from below by appropriate functions of $k$ and $\Delta(G)$, the corresponding thresholds admit the following bounds: $\ra_k'(G) \ge k$,  $\rb_{k-1}'(G) \leq k$, $\ra_k(G) \ge k$ and $\rb_{k-1}(G) \leq k$. In Section~\ref{sec:upper-bound-a1}, we introduce a new graph invariant, the local domination number, $\widetilde{\gamma}_\ell(G)$, of a graph $G$. We prove that $\widetilde{\gamma}_\ell(G)$ is an upper bound for $\ra'_\ell(G)$ and $\ra_\ell(G)$ in any graph $G$ with $\delta(G)\ge \ell$, which in turn implies that $\ra_1'(G)\le k-1$ holds for any induced $K_{1,k}$-free graph $G$ with no isolated vertices. In Section~\ref{sec:line-graph}, we prove upper bounds on Dominator's thresholds and lower bounds on Staller's thresholds in line graphs. The result is also instrumental in providing the sharpness of the bounds from Section~\ref{sec:large-order} as it enables us to show that for sufficiently large order there are line graphs $G$ for which $\ra_k(G)=\ra_k'(G) = k$ and  $\rb_{k-1}(G) =\rb_{k-1}'(G) = k$. In Section~\ref{sec:grids}, we prove several exact values of the thresholds in Cartesian grids, both finite and infinite. Finally, in Section~\ref{sec:star-partition}, we introduce the concept of star-partition width, which is closely related to a $[1,k]$-factor of a graph that goes back to~\cite{Lovasz-1970}; see also~\cite{heinrich-1990}. We prove that the star partition width of a nontrivial graph $G$ is an upper bound on $\ra_1'(G)$, while in trees this bound is attained. Using a different language the result about trees was obtained in~\cite{bagdas-2024}, while our proof is completely different and somewhat shorter.

\paragraph{Standard definitions.} A {\em leaf} of a graph $G$ is a vertex of degree $1$.  
For $k \ge 1$, the {\em $k$-star} is the complete bipartite graph $K_{1,k}$. If $u$ and $v$ are two vertices of a graph, ${\rm dist}(u,v)$ denotes their distance. As usual, $n(G)$, $\delta(G)$, and $\Delta(G)$ respectively stand for the \emph{number of vertices, minimum and maximum vertex degree} in a graph $G$, while $\nu(G)$,  $\tau(G)$, $\alpha(G)$, and $\gamma(G)$ denote the \emph{matching number, vertex cover number, independence number and domination number} of $G$. If $S\subseteq V(G)$, then $G[S]$ denotes the subgraph of the graph $G$ induced by $S$. For a vertex $v \in V(G)$, the \emph{open and closed neighborhoods} of $v$ are denoted by $N(v)$ and $N[v]$, respectively.
If $k$ is a positive integer, then $[k]$ stands for the set $\{1,\dots , k\}$. 

The \emph{line graph} $L(G)$ of a graph $G$ is defined as follows. The vertices of $L(G)$ represent the edges of $G$, and two vertices are adjacent in $L(G)$ if the corresponding edges share a vertex in $G$. A graph $F$ is called line graph if there is a graph $G$ such that $F=L(G)$. 

 Let $G$ and $F$ be two graphs. The \emph{Cartesian product} $G \Box F$ is the graph on the vertex set $V(G) \times V(F)$ where two vertices $(x,y)$ and $(x',y')$ are adjacent if either $x=x'$ and $yy' \in E(F)$, or $y=y'$ and $xx' \in E(G)$. For a vertex $x \in V(G)$, we denote by $^x\!F$ the subgraph of $G \Box F$ induced by the vertex set $\{ (x,y): y \in V(F)\} $. For every $x \in V(G)$, this subgraph $^x\!F$ is isomorphic to $F$ and called an \emph{$F$-layer} in $G \Box F$. A \emph{$G$-layer} in the product is defined analogously and denoted by $G^y$ for a vertex $y \in V(F)$.

%%%%%%%%%%%%%%%%%%%%%%%%%%%%%%%%%%%%%%
\section{Preliminaries} 
\label{sec:preliminaries}

In this section, we make some basic observations and prove simple but sharp upper bounds on the thresholds considered. Note that in Maker-Breaker positional games it is never in favor of a player to pass a move. Hence, the following observation is not surprising and is easy to see. 
\begin{observation}
\label{obs:basic}
Let $G$ be a graph and $a$ and $b$ be positive integers. \begin{enumerate}[(i)]
    \item If $W(G,a,b)=S$, then $W(G,a,b+1)=S$; if $W'(G,a,b)=S$, then $W'(G,a,b+1)=S$.
     \item If $W(G,a,b)=D$, then $W(G,a+1,b)=D$; if $W'(G,a,b)=D$, then $W'(G,a+1,b)=D$.
    \item If $W(G,a,b)=S$, then $W'(G,a,b)=S$.
    \item If $W'(G,a,b)=D$, then $W(G,a,b)=D$.
\end{enumerate}
\end{observation}

Due to the monotonicity of the functions $W$ and $W'$, as noted in Observation~\ref{obs:basic}, items (i) and (ii), the following definitions arise. {\em Dominator's threshold}, ${\rm a}_\ell(G)$, in the Dominator-start biased MBD game played on $G$ is 
defined as $\min\{a:  W(G,a,\ell)=D\}$, while {\em Dominator's threshold} in the Staller-start  biased MBD game is denoted by ${\rm a}_\ell'(G)$  and is defined as $\min\{a:  W'(G,a,\ell)=D\}$. Similarly, {\em Staller's thresholds} in both versions of the game are defined by
${\rm b}_\ell(G)=\min\{b:  W(G,\ell,b)=S\}$ and
${\rm b}_\ell' (G)=\min\{b:  W'(G,\ell,b)=S\}$.
If the minimum in the definition does not exist, we set $\infty$ as the value of the threshold. For example, $\ra_2'(T)=\infty$ for every tree $T$ as Staller wins the $(\ell, 2)$-biased S-game for every $\ell$ by playing a leaf and its neighbor in her first move.
By Observation~\ref{obs:basic}, items (iii) and (iv), we infer that 
$${\rm b}_{\ell}(G)\ge {\rm b}_\ell'(G) \textrm{  and }{\rm a}_\ell'(G)\ge {\rm a}_\ell(G)$$ hold for any graph $G$. 

By definitions and Observation~\ref{obs:basic} (i) and (ii), we obtain the following statements. 
\begin{proposition} \label{prop:W}
The following equivalences are true for every graph $G$ and positive integers $i$ and $j$.
\begin{enumerate}[(i)]
    \item $W(G,i,j)=D$, if and only if $\ra_j(G) \leq i$, and if and only if $\rb_i(G)\ge j+1$.
    \item $W'(G,i,j)=D$, if and only if $\ra_j'(G) \leq i$, and if and only if $\rb_i'(G)\ge j+1$.
    \item $W(G,i,j)=S$, if and only if $\ra_j(G) \ge i+1$, and if and only if $\rb_i(G)\le j$.
    \item $W'(G,i,j)=S$, if and only if $\ra_j'(G) \ge i+1$, and if and only if $\rb_i'(G)\le j$.
\end{enumerate}    
\end{proposition}

\medskip
Observation~\ref{obs:basic} and Proposition~\ref{prop:W} together imply that $\ra_i(G) \ge \ra_j(G)$, $\ra_i'(G) \ge \ra'_j(G)$, $\rb_i(G) \ge \rb_j(G)$, and $\rb_i'(G) \ge \rb'_j(G)$  hold when $i >j$. We next show that there exist trivial sharp upper bounds for all four thresholds.
\begin{proposition} \label{prop:trivial-bounds}
     The following inequalities hold for every graph $G$ and positive integers $a$, $b$.
    \begin{enumerate}[(i)]
        \item ${\rm b}'_a(G) \leq \delta(G)+1$;
        \item ${\rm b}_a(G) \leq \Delta(G)+1$ when $a < \gamma(G)$;
        \item ${\rm a}_b(G)\leq \gamma(G)$;
        \item ${\rm a}'_b(G) \leq b \cdot\Delta(G)$ when $b \le \delta(G)$.
    \end{enumerate}
    Moreover, all bounds are tight.
\end{proposition}
\begin{proof}
   (i) Consider the $(a, \delta(G)+1)$-biased S-game on $G$. Staller starts by playing a vertex of minimum degree and all its neighbors, and thus she wins the game.  Therefore, ${\rm b}'_a(G) \leq \delta(G)+1$ for every positive integer $a$.

   To show the sharpness of the upper bound, we first introduce a graph $G_{n,k}$, for every $n\ge 2$ and $k \ge 3$. This graph is obtained from $n$ copies of $K_k$ by inserting $n-1$ edges between the components to make the graph connected such that $\Delta(G_{n,k})= k$ and every copy of $K_k$ is incident to at most two such edges, see Fig.~\ref{Gnk}. Let $G=G_{a,k}$ with $a\ge 2$ and $k\ge 3$. The latter condition ensures that $\delta(G)=k-1$ holds. 
   By the proof above, ${\rm b}'_a(G) \le k$.
   We now show that $W'(G,a,k-1)=D$. Since $\delta(G)=k-1$, Staller cannot claim a closed neighborhood in $G$ in her first move, and Dominator can win the game in the next move by selecting one vertex from each copy of $K_k$.
   Thus ${\rm b}'_a(G) = k=\delta(G)+1$.

   We remark that, alternatively, any graph $G$ with $\gamma(G)= a $ and with $\delta(G) +1$ pairwise disjoint minimum dominating sets is a sharp example for (i) as ${\rm b}'_a(G) =\delta(G)+1$.

\begin{figure}
\begin{center}
	\begin{tikzpicture}[]
        %edges	
	\draw (0,0) -- (11,0);
	\draw (0,0)--(1,1.5)--(2,0);
        \draw (3,0)--(4,1.5)--(5,0);
        \draw (6,0)--(7,1.5)--(8,0);
        \draw (9,0)--(10,1.5)--(11,0);
	
	%vertices

	\draw [fill=black] (1,1.5) circle (5pt);
	\draw [fill=black] (0, 0) circle (5pt);
	\draw [fill=black] (2, 0) circle (5pt);
	\draw [fill=black] (3, 0) circle (5pt);
	\draw [fill=black] (5, 0) circle (5pt);
	\draw [fill=black] (6, 0) circle (5pt);
	\draw [fill=black] (8, 0) circle (5pt);
	\draw [fill=black] (9, 0) circle (5pt);
        \draw [fill=black] (11, 0) circle (5pt);
       \draw [fill=black] (4,1.5) circle (5pt);
       \draw [fill=black] (7,1.5) circle (5pt);
       \draw [fill=black] (10,1.5) circle (5pt);

	\end{tikzpicture} \\
	\caption{The graph $G_{4,3}$.}\label{Gnk}
	
\end{center}
\end{figure}

   (ii) Now consider the $(a, \Delta(G)+1)$-biased D-game on $G$ and assume that  $a < \gamma(G)$. After Dominator's first move, there is an undominated vertex $u$. As $|N[u]| \leq \Delta(G) +1 $, Staller can play all vertices from $N[u]$ and wins the game with her first move. We conclude that ${\rm b}_a(G) \leq \Delta(G)+1$. 
   
   To show that the upper bound is sharp, we define $F_{a,n}=C_{a+1} \square K_n$ for $a \ge n \ge 5$. Then $\Delta(F_{a,n}) = n+1$ and $\gamma(F_{a,n})=a+1$. 
   Since ${\rm b}_a(G) \leq n+2$, it remains to show that Dominator has a winning strategy in the $(a,n+1)$-biased D-game on $F_{a,n}$. Dominator selects one vertex from each of $a$ arbitrary $K_n$-layers in his first move, and then only one $K_n$-layer, say $^{v_1}\! K_n$, remains unplayed. Further, at most $n-1$ vertices, all from $^{v_1}\! K_n$, remain undominated after Dominator's first move. Since every closed neighborhood contains $n+2$ vertices, Staller cannot claim a closed neighborhood of a vertex in her next move. Then, Dominator can win the game in his second move by choosing an unplayed neighbor for each undominated vertex from $^{v_1}\! K_n$.
   Therefore, $W(F_{a,n}, a, n+1)=D$ that implies  ${\rm b}_a(F_{a,n}) = n+2=\Delta(F_{a,n})+1$.

   (iii) In the D-game, if Dominator can play $\gamma(G)$ vertices per move, then he can form a dominating set and win the game in his first move. Thus ${\rm a}_b(G)\leq \gamma(G)$.

   Next, we prove the sharpness of the upper bound. Consider again the graph $G=G_{n,b}$ defined in the proof of part (i). Assuming that $n\ge 2$ and $b \ge 3$, we have $\Delta(G)=b$, $\delta(G)=b-1$, and $\gamma(G)=n$. 
  We show that $W(G,n-1,b)=S$. After Dominator's first move, there is an unplayed $K_b$-subgraph.  In her first move, Staller plays all vertices from this clique. As there is a vertex $u$ of degree $b-1$ in this clique, Staller claims the entire neighborhood $N[u]$ with this move. Thus, Staller wins the $(n-1,b)$-biased D-game, and we infer that ${\rm a}_b(G) = n=\gamma(G)$.

    Alternatively, any $(b-1)$-regular graph $G$ is a sharp example  as ${\rm a}_b(G) = \gamma(G)$ can be proved analogously to the case of $G=G_{n,b}$.

   (iv) Assume that $b \le \delta(G)$. We consider the $( b \cdot\Delta(G), b)$-biased S-game on $G$. Since  $b \le \delta(G)$, Staller cannot claim a closed neighborhood of a vertex in her first move. Then Dominator replies by selecting all unplayed neighbors of the $b$ vertices played by Staller (and some further unplayed vertices if needed). If Dominator continues applying this strategy, Staller cannot claim a closed neighborhood in $G$. Hence, Dominator wins the game and we conclude that ${\rm a}'_b(G) \leq b \cdot\Delta(G)$.

   As sharp examples with $b=1$, we refer to the odd paths $P_n$. It is known that Staller wins the $(1,1)$-biased S-game on $P_n$; see~\cite{duchene-2020}. Thus we have ${\rm a}'_1(P_n) = 2= 1 \cdot\Delta(P_n)$ if $n$ is odd. 
\end{proof}
\begin{remark}
    Proposition~\ref{prop:W} and parts (i)-(iii) of Proposition~\ref{prop:trivial-bounds} imply the following statements:
    \begin{enumerate}[(i)]
        \item ${\rm a}'_{\delta(G)+1}(G) =\infty$;
        \item ${\rm a}_{\Delta(G)+1}(G) =\gamma(G)$;
        \item ${\rm b}_{\gamma(G)}(G)=\infty$.
        \end{enumerate}
\end{remark}

%%%%%%%%%%%%%%%%%%%%%%%%%%%%%%
\section{Graphs of large order}
\label{sec:large-order}

In this section, we show that Staller can win the $(a,b)$-biased MBD game in $G$ if $a<b$ and the order of the graph is large enough. Intuitively, a large enough order enables Staller to initially choose a large set of pairwise disjoint closed neighborhoods for which she can ensure, since she selects more vertices in each move than Dominator, that enough of them will not be played by Dominator. At the end of the game this will enable Staller to select all vertices of one of the mentioned closed neighborhood.

\begin{theorem} \label{thm:large}
    Let $k \ge 2$ be an integer and $G$ be a graph of maximum degree $\Delta(G)= \Delta$.
    \begin{enumerate}[(a)]
        \item If $n(G) \ge k^{\Delta-k+1}(\Delta^2+1)$, then $\ra_k'(G) \ge k$ and $\rb_{k-1}'(G) \leq k$.
         \item If $\Delta \ge k$ and $n(G) \ge (k^{\Delta-k+1}+1)(\Delta^2+1)$, then $\ra_k(G) \ge k$ and $\rb_{k-1}(G) \leq k$.
         \item If $\Delta < k$ and $n(G) > (k-1)(\Delta+1)$, then $\ra_k(G) \ge k$ and $\rb_{k-1}(G) \leq k$.
    \end{enumerate}
\end{theorem}
\begin{proof}
We consider the $(k-1,k)$-biased S-game and D-game on $G$ and prove that Staller wins under the conditions given on $n(G)$.
\medskip

(a) 
Assume that $n(G) \ge k^{\Delta-k+1}(\Delta^2+1)$.
If $\Delta < k$, then Staller wins the S-game on $G$ by claiming the closed neighborhood of a vertex in her first move. From now on, we may suppose that $\Delta \ge k$. To describe Staller's strategy to win, we split the game into $\Delta-k+2$ parts. For every $i \in [\Delta-k+1]$, Part ($i$) finishes with Dominator's move after $k^{\Delta-k+1-i}$ (full) moves of Dominator and Staller  (except when Staller wins during this part). In each Part (i) of the game Staller has a special goal, which we denote by S($i$). The goal of the first part is as follows.

\paragraph{S($1$)} Staller's first goal is to finish Part ($1$) with $r_1 \ge k^{\Delta-k}$ vertices $s^1_1, \dots, s^1_{r_1}$ each played by Staller such that the neighborhoods $N[s^1_1], \dots ,N[s^1_{r_1}] $ are pairwise disjoint and Dominator played none of the vertices from $\bigcup_{j=1}^{r_1} N[s^1_j]$.
\medskip

For every two vertices $u$ and $v$ in $G$, the neighborhoods $N[u]$ and $N[v]$ intersect only if ${\rm  dist}(u,v) \leq 2$. For a vertex $u \in V(G)$, at most $\Delta^2+1$ such vertices exist. As $n(G)\ge k^{\Delta-k+1}(\Delta^2+1)$, Staller can choose a set $A$ of  $k^{\Delta-k+1}$ vertices such that the neighborhoods of the vertices are pairwise disjoint. During her first $k^{\Delta-k}$ moves in Part (1), Staller plays the $k\cdot k^{\Delta-k}$ vertices in $A$. (If Dominator plays a vertex from $A$, then Staller plays an arbitrary vertex instead of it.)  During his first $k^{\Delta-k}$ moves, Dominator plays $(k-1)k^{\Delta-k}$ vertices. This way, there remain at least $k^{\Delta-k}$ vertices in $A$ that were played by Staller and their neighborhoods contain no vertex played by Dominator. These vertices satisfy Staller's aim S($1$).
\medskip

In the continuation of the game, Staller's goal in Part ($i$) is the following.
\paragraph{S($i$)} Staller aims to finish Part ($i$) with either winning during this part or having $r_i \ge k^{\Delta-k+1-i}$ vertices $s^i_1, \dots, s^i_{r_i}$ such that $N[s^i_1], \dots, N[s^i_{r_i}]$ are pairwise disjoint, Staller played $i$ vertices from each closed neighborhood, and Dominator played no vertex from  $\bigcup_{j=1}^{r_i} N[s^i_j]$. (Note that $\{s^i_1, \dots, s^i_{r_i} \} \subseteq \{s^{i-1}_1, \dots, s^{i-1}_{r_{i-1}} \}$.) 
\medskip

By induction on $i$, we prove that Staller can achieve her goal for every $i \in [\Delta-k+1]$. We have already seen that it is possible for $i=1$. Suppose that $2 \le i \le \Delta-k+1$ and Staller does not win in Part ($i$). By our hypothesis, at the end of Part ($i-1$) there are $r_{i-1} \ge r_{i-1}'= k^{\Delta-k+2-i}$ vertices that satisfy the properties in S($i-1$). Among these, arbitrarily select vertices $s^{i-1}_1, \dots, s^{i-1}_{r_{i-1}'}$.  Then, during her $k^{\Delta-k+1-i}$ moves in Part ($i$), Staller plays one vertex from each neighborhood $N[s^{i-1}_1], \dots , N[s^{i-1}_{r_{i-1}'}]$ (if there is an unplayed vertex). In the meantime, Dominator plays only  $(k-1)k^{\Delta-k+1-i}$ vertices that leave at least $k^{\Delta-k+1-i}$ disjoint neighborhoods unplayed by Dominator such that Staller played $i$ vertices from each of them. It shows that Staller can achieve her goal.
\medskip

At the end of Part ($\Delta-k+1$), if Staller did not win earlier, S($\Delta-k+1$) is satisfied and we have $r_{\Delta-k+1} \ge 1$ vertices with the prescribed property. For such a vertex $x=s^{\Delta-k+1}_1$, the neighborhood $N[x]$ already contains  $\Delta-k+1$ vertices played by Staller and no vertex from $N[x]$ is played by Dominator. Therefore, Staller can win the game by choosing the remaining $|N[x]|- (\Delta-k+1) \le k$ vertices from $N[x]$ in her next move.  

It proves $W'(G,k-1,k)=S$ that in turn implies $\ra_k'(G) \ge k$ and $\rb_{k-1}'(G) \leq k$   by Proposition~\ref{prop:W} (iv).
\medskip

(b) 
If $\Delta \ge k$ and $n(G) \ge (k^{\Delta-k+1}+1)(\Delta^2+1)$, the proof for the D-game is analogous to the proof of part $(a)$. We remark that if Dominator's first move is playing $d_1, \dots d_{k-1}$, then Staller may consider the graph  $G'=G - \bigcup_{j=1}^{k-1} N[d_j]$ and apply the strategy from the proof of part (a) for $G'$. As 
$$n(G') \ge n(G) - (k-1)(\Delta +1) \ge n(G)-(\Delta^2-1) > k^{\Delta-k+1}(\Delta^2+1),
$$ 
Staller can win the game, and the statement follows by Proposition~\ref{prop:W} (iii).
\medskip

(c)  As $n(G) > (k-1)(\Delta +1)$, the $(k-1)$ vertices played by Dominator in his first move are not enough to dominate $G$. On the other hand,  $\Delta < k$ ensures that Staller can choose an undominated vertex $v$ after Dominator's move and can play all vertices from $N[v]$. This way, Staller wins and $W(G,k-1,k)=S$ implies the statement.
\end{proof}

\begin{remark}\label{rem:large}
    We note that Staller's winning strategy in the proof of Theorem~\ref{thm:large} also works if $G$ is an infinite graph with a finite maximum degree $\Delta$. Hence, the lower and upper bounds on the thresholds stated in Theorem~\ref{thm:large} remain true for infinite graphs with finite maximum degrees. 
\end{remark}
\begin{remark}
    Corollary~\ref{cor:line-graph-large} in Section~\ref{sec:line-graph} will show that Theorem~\ref{thm:large} is sharp for infinitely many graphs. 
\end{remark}

%%%%%%%%%%%%%%%%%%%%%%%%%%%%%%%%%%%%%%%%%%%%%%%%%%%%%%%%%%%%%%%%%%%%%%%%%%%%%%%
\section{Local domination and an upper bound on $\ra_\ell'(G)$}
\label{sec:upper-bound-a1}

In this section, we study the S-game where Staller's bias equals $\ell$, and prove an upper bound on Dominator's threshold ${\rm a}_\ell'(G)$. Since $\ra_\ell(G) \le \ra_\ell'(G)$ holds for every graph $G$, the established upper bound is also valid for $\ra_\ell(G)$. 

We introduce a parameter called the $\ell$-local domination number of a graph $G$. Let $G$ be a graph with $\delta(G) \geq \ell$. For $S \subseteq V(G)$ with $|S|=\ell$, 
let $\widetilde{\gamma}_\ell(G,S)$ be the minimum number of vertices in $V(G) \setminus S$ that can dominate $N(S)\cup \{v \in S:\, {\rm{deg}}(v)=\ell\}$, where $N(S)=\bigcup_{v \in S}N(v) $. Then, the \emph{$\ell$-local domination number}, $\widetilde{\gamma}_\ell(G)$, of $G$ is defined as
 $$ \widetilde{\gamma}_\ell(G)= \max_{S \subseteq V(G),|S|=\ell} \{\widetilde{\gamma}_\ell(G,S) \}. $$

To illustrate the definition, note that $\widetilde{\gamma}_1(C_n)=2$ when $n\ge 5$ and $\widetilde{\gamma}_2(C_n)=4$ when $n\ge 10$. As an additional example, consider the graph $G$ depicted on Fig.~\ref{fig:localDom}. If $S =\{z\}$ or $S=\{y_1\}$ or $S=\{y_2\}$, then $\widetilde{\gamma}_1(G,S)=2$. If $S=\{x_i\}$ for some $i \in [4]$, then $\widetilde{\gamma}_1(G,S)=1$. Hence $\widetilde{\gamma}_1(G)=2.$

\begin{figure}[ht!]
\begin{center}
	\begin{tikzpicture}[scale=1.1]
        %edges	
	\draw (0,0) -- (2,0);
	\draw (0,2)--(0,0)--(-2,2);
        \draw (2,2)--(2,0)--(4,2);
        \draw (-2,2)--(1,4)--(0,2);
        \draw (2,2)--(1,4)--(4,2);
	
	%vertices

	\draw [fill=black] (0, 0) circle (5pt);
	\draw [fill=black] (2, 0) circle (5pt);
	\draw [fill=black] (-2, 2) circle (5pt);
	\draw [fill=black] (0, 2) circle (5pt);
	\draw [fill=black] (2, 2) circle (5pt);
	\draw [fill=black] (4, 2) circle (5pt);
	\draw [fill=black] (1, 4) circle (5pt);

	\draw (1,4.4) node{$z$}; 
    \draw (-2.4,2.2) node{$x_1$};
    \draw (-0.4,2.2) node{$x_2$};
    \draw (1.6,2.2) node{$x_3$};
    \draw (4.3,2.2) node{$x_4$};
    \draw (-0.4,-0.4) node{$y_1$};
    \draw (2.4,-0.4) node{$y_2$};
	
	\end{tikzpicture} \\
	\caption{A graph $G$.}\label{fig:localDom}
\end{center}
\end{figure}

 \begin{theorem} \label{thm:k-local-dom-a1'}
     If $G$ is a graph with minimum degree at least $\ell$, then ${\rm a}_\ell'(G) \leq \widetilde{\gamma}_\ell(G)$.
 \end{theorem}  
 
 \begin{proof}
  It suffices to prove that Dominator has a strategy to win the $(\widetilde{\gamma}_\ell(G), \ell)$-biased S-game on $G$. Let us consider the following strategy for Dominator to follow:
 \begin{itemize}
 \item[$(\star)$] After each move of Staller selecting vertices from $S$, where $S\subset V(G), |S|=\ell$, Dominator makes sure that every vertex in $N(S)\cup \{v \in S:\, {\rm{deg}}(v)=\ell\}$ is dominated after Dominator's next move. 
 \end{itemize}
 As $\delta(G) \geq \ell$, Staller cannot win the game with her first move $S_1=\{s^1_1,\ldots , s^1_\ell\}$. According to strategy $(\star)$, Dominator replies by playing $\widetilde{\gamma}_\ell(G)$ vertices that dominate $N(S_1)\cup \{v \in S_1:\, {\rm{deg}}(v)=\ell\}$. It is possible as $\widetilde{\gamma}_\ell(G,S_1) \leq \widetilde{\gamma}_\ell(G) $. (In case of strict inequality, Dominator selects some additional vertices arbitrarily.) Now, we proceed by induction.

Suppose that Dominator could play according to strategy $(\star)$ during his first $i$ moves, and the game is not yet over. Via three claims, we prove that Staller is not able to win with her $(i+1)$st move and that Dominator can choose his $(i+1)$st move complying with $(\star)$.
 \medskip
 
\noindent \textbf{Claim 8.A. }
 \textit{Assume that no vertex from  $N[v]$ was played by Dominator during his first $i$ moves. If $\deg(v)=\ell$, then all vertices in $N[v]$ are unplayed, and otherwise all vertices from the open neighborhood $N(v)$ are unplayed after Dominator's $i$th move.}
 
\smallskip 

 \noindent\emph{Proof.}
 First, let deg$(v)=\ell$ and suppose for a contradiction that Staller has already played a vertex from $N[v]$. In this case, let $u$ be the first vertex from $N[v]$ played by Staller, say in Staller's $j$th move, where $j\le i$. If $u=v$, then, according to the hypothesis and since deg$(v)=\ell$, Dominator replied by playing a neighbor $v'$ of $v$ (and maybe some further vertices) as his $j$th move.
 If $u\neq v$, then Dominator dominated $v \in N(u)$ after Staller's move $u$. In either case, we get a contradiction with the assumption that Dominator has not played any vertex in $N[v]$ during his first $i$ moves.  

If deg$(v) > \ell$, suppose that Staller played earlier a neighbor $u$ of vertex $v$. As Dominator followed strategy $(\star)$, in the next turn, he dominated the entire $N(u)$. It included dominating $v$ itself that is playing a vertex from $N[v]$. This again contradicts the assumption in Claim 8.A and finishes the proof of the statement. \hfill$(\square)$  
 \medskip
  
\noindent \textbf{Claim 8.B. }
 \textit{Staller cannot win the game with her $(i+1)$st move.}
 
\smallskip 

 \noindent\emph{Proof.}
 Suppose to the contrary that Staller wins in her $(i+1)$st move $S_{i+1}=\{s_1^{i+1},\ldots ,s_\ell^{i+1}\}$ by achieving that all vertices in $N[u]$ are played by her. Since $|N[u]| > \ell$, at least one vertex $v \in N[u]$ was chosen by Staller during her first $i$ moves.
Since no vertex from $N[u]$ was played by Dominator, it follows by Claim 8.A that deg$(u)>\ell$ and $u$ is the only vertex in $N[u]$ played by Staller in her first $i$ moves. This is a contradiction since Staller cannot claim the whole $N(u)$ in his $(i+1)$st move because $|N(u)|>\ell$. 
 \hfill$(\square)$  
 \medskip

 \noindent \textbf{Claim 8.C. }
 \textit{Dominator can play according to strategy $(\star)$ in his $(i+1)$st move.}  
\smallskip 

 \noindent\emph{Proof.}
 Let $S_{i+1}=\{s_1^{i+1},\ldots ,s_\ell^{i+1}\}$ be Staller's $(i+1)$st move. Let $R \subseteq V(G) \setminus S_{i+1}$ be a set with $|R| \leq \widetilde{\gamma}_\ell(G)$ such that $R$ dominates the entire $N(S_{i+1}) \cup \{v \in S_{i+1}:\, {\rm{deg}}(v)=\ell\}$. If all vertices in $R$ are unplayed, then Dominator plays these vertices in his $(i+1)$st move. If a vertex $x \in R$ was played earlier by Dominator, then he does not have to play it in the $(i+1)$st move (he may replace it with an arbitrary vertex). If a vertex $y \in R$ was played by Staller in her $j$th move, where $j < i+1$, we have two cases. 
\begin{itemize}
    \item 
 First, suppose that $y \notin N(S_{i+1})$. Since $j < i+1$, the induction hypothesis implies that Dominator dominated all neighbors of $y$ with his $j$th move. Hence, in his $(i+1)$st move, vertex $y$ is not needed in $R$ to dominate $N(y)$. 
\item  In the other case, let $y \in N(S_{i+1})$. Then all vertices in $N(y)$ are dominated by Dominator after his $j$th move due to $(\star)$. However, it might happen that $y$ itself is not dominated by the earlier moves of Dominator. By Claim 8.A, either Dominator has already played a vertex from $N[y]$ and then $y$ is dominated, or all vertices in $N(y)$ were unplayed before the move $S_{i+1}$ and deg$(y)> \ell$. In this case, Dominator can choose an arbitrary vertex from $N(y) \setminus S_{i+1}$ to play instead of $y$ in his $(i+1)$st move. Note that only the vertices in $S_{i+1} \cap N(y)$ are played from $N(y)$ and, as $|N(y)| \ge \ell+1$, Dominator can choose an unplayed vertex $z \in N(y)\setminus S_{i+1}$ and plays $z$ instead of $y$ in his $(i+1)$st move.   \hfill$(\square)$  
\end{itemize} 
  Claims 8.B and 8.C prove that Dominator can follow strategy~$(\star)$ throughout the $(\widetilde{\gamma}_\ell(G), \ell)$-biased S-game on $G$, and that Staller is never able to claim a closed neighborhood to win the game.
This implies ${\rm a}_\ell'(G) \leq \widetilde{\gamma}_\ell(G)$ as stated in the theorem.  
 \end{proof}

The $1$-local domination number of a graph $G$ is also called \emph{local domination number} and denoted by $\widetilde{\gamma}(G)$. If $v$ is a leaf in $G$, then both $N(v)$ and $N[v]$ can be dominated with only one vertex from $V(G) \setminus \{v\}$. Therefore, the definition of $\widetilde{\gamma}(G)$ can be simplified as follows. Let $G$ be a graph without isolated vertices. For a vertex $v \in V(G)$, let $\widetilde{\gamma}(G,v)$ be the minimum number of vertices that can dominate $N(v)$ in the graph $G-v$. Then, the local domination number is defined as
 $$ \widetilde{\gamma}(G)= \max_{v \in V(G)} \{\widetilde{\gamma}(G,v) \}. $$ 

Theorem~\ref{thm:k-local-dom-a1'} with $\ell=1$ gives the upper bound $\ra'_1(G) \leq \widetilde{\gamma}(G)$ for every graph with $\delta(G) \ge 1$. 
For a claw-free graph $G$ (i.e., when no induced subgraph of $G$ is a star $K_{1,3}$), it holds by definition that $\alpha(G[N(v)]) \leq 2$ for every $v \in V(G)$. Consequently, the local domination number of a claw-free graph $G$ is at most $2$. Similarly, if $G$ does not contain any induced star subgraph $K_{1,k}$, then $\widetilde{\gamma}(G) \leq k-1$. Thus Theorem~\ref{thm:k-local-dom-a1'} directly implies the following statements. 
\begin{corollary} \label{cor:claw-free}
    If $G$ is an induced $K_{1,k}$-free graph with $\delta(G) \ge 1$, then ${\rm a}_1'(G) \leq k-1$ holds. In particular, if $G$ is a claw-free graph, then ${\rm a}_1'(G) \leq 2$.
\end{corollary}
For claw-free graphs, the upper bound ${\rm a}_1'(G) \leq 2$ is sharp. For example, it was shown in~\cite{duchene-2020} that $W'(P_{2k+1},1,1)=S $ holds for every positive integer $k$ that implies ${\rm a}_1'(P_{2k+1}) = 2$ by Corollary~\ref{cor:claw-free}. Another claw-free graph $P_{2k+1}^+$ can be obtained from the odd path $P_{2k+1}: v_1\dots v_{2k+1}$ if we add the edges between $v_{2i}$ and $v_{2i+2}$ for every $i \in [k-1]$. The strategy described in the proof of~\cite[Theorem 5.2]{bujtas-2021} for Staller's speedy win on $P_{2k+1}$ can be applied to the extended graph $P_{2k+1}^+$. It shows that $W'(P_{2k+1}^+,1,1)=S $, and then ${\rm a}_1'(P_{2k+1}^+) = 2$ holds for every positive integer $k$.
\medskip

%%%%%%%%%%%%%%%%%%%%%%%%%%%%%%%%%%%%%%%%%%%%%%%%%%%%%%%%%%%%%%%%%%%%%%
\section{Thresholds for line graphs} \label{sec:line-graph}

As every line graph is a claw-free graph, Corollary~\ref{cor:claw-free} implies ${\rm a}_1'(G) \leq 2$ for each line graph $G$. We show that, in most of the cases, the exact value for ${\rm a}_1'(G)$ is $1$. Further, we prove that $\ra_k(G) = \ra_k'(G)= k$ and $\rb_k(G)= \rb_k'(G) = k+1$ holds for a large subclass of line graphs.

  Consider a collection $\cF=\{F_1, \dots , F_n\}$ of sets. A \emph{set of distinct $t$-representatives}  (SDR$^t$) for $\cF$ is a set $\{r_i^j: i \in [n], \enskip j \in [t]\}$ of $nt$ different elements such that $r_i^j \in F_i$ holds for every pair $(i,j)$ with $i \in [n]$ and $ j \in [t]$. Note that SDR$^1$ gives exactly the well known concept of the set of distinct representatives. For larger values of $t$, SDR$^t$ is obtained by selecting $t$ elements from each set in $\cF$ so that all selected elements are distinct. 
  Using Hall's celebrated theorem from~\cite{hall-1935}, the following generalized version holds; see \cite{Mirsky-1967}. 
  \begin{theorem} \label{thm:t-hall}
      A family $\cF=\{F_1, \dots , F_n\}$ has a  set of distinct $t$-representatives if and only if, for every subfamily $\cF' \subseteq \cF$, the union $\bigcup_{F_i \in \cF'} F_i  $ contains at least $t\,|\cF'|$ elements.  
  \end{theorem}

\begin{theorem} \label{thm:line-graph-general}
    If $G$ is the line graph of a graph $H$ with $\delta(H) \ge 2t$ and $k \le 2t-1$, then $\ra_k(G)\leq \ra_k'(G)\le k$ and $\rb_k(G)\geq \rb_k'(G) \geq k+1$.  
\end{theorem}
\begin{proof}
 It is enough to prove that, under the conditions in the theorem, Dominator wins in the $(k,k)$-biased S-game on $G=L(H)$. In the proof, we consider a set of cliques $\cS=\{Q_1, \dots , Q_n\}$ in $G$ such that $n=n(H)$ and, for every $i \in [n]$, clique $Q_i$ corresponds to the set of edges incident with vertex $u_i \in V(H)$. By the condition $\delta(H) \ge 2t$, we have $|Q_i| \ge 2t$ for every $Q_i \in \cS$. By definition, every vertex $v \in V(G)$ is contained in exactly two cliques from $\cS$. 
  \medskip

  We first show that there is a set of distinct $t$-representatives for the set system $\cS$. By Theorem~\ref{thm:t-hall}, it is enough to prove that $|\bigcup_{Q \in \cS'} Q| \ge t\,|\cS'|$ holds for every subfamily $\cS' \subseteq \cS$. Consider the sum $x = \sum_{Q \in \cS'} |Q| $. As every clique $Q \in \cS$ contains at least $2t$ vertices, we have $x \ge 2t\,|\cS'|$. On the other hand, every vertex from $\bigcup_{Q \in \cS'} Q$ is covered by at most two cliques from $\cS'$ and hence, $2 |\bigcup_{Q \in \cS'} Q| \ge x $ holds. This implies $|\bigcup_{Q \in \cS'} Q| \ge t|\cS'|$  for every $\cS' \subseteq \cS$, and we may conclude that $\cS$ admits an SDR$^t$.
  \medskip

  Fixing an SDR$^t$ in $\cS$, for every $i \in [n]$, let $R_i=\{r_i^1, \dots , r_i^t\}$ be the set of the $t$ representative vertices assigned to the clique $Q_i$. Consider the following strategy of Dominator in the $(k,k)$-biased S-game on $G$. If Staller played vertices $x_1, \dots, x_k$ in a move, then Dominator replies by choosing vertices $y_1, \dots, y_k$ (one by one) such that the following is true for every $i \in [k]$.
  \begin{itemize}
      \item If $x_i \in R_j$, and there is an unplayed vertex $y_i\in Q_j$, then Domination plays $y_i$. If there no such vertex in $Q_j$, Dominator selects a neighbor of $x_i$ if possible. Finally, if all neighbors of $x_i$ have already been played, then he selects a vertex randomly. 
      \item If $x_i$ is not a representative vertex for any clique $Q_j$, Dominator plays a neighbor $y_i$ of $x_i$. If all neighbors have been played before, Dominator chooses an arbitrary vertex. 
  \end{itemize}
  Assuming that Dominator follows this strategy, we may observe the following: 
  \paragraph{Claim 11.A.}
      \textit{If Staller plays a  vertex $x_i$ in a move, then, after this move, there is a vertex $z \in N_G[x_i]$ that is either unplayed or  played by Dominator.}\\
      \textit{Proof.} Assume to the contrary that after Staller's move, every vertex in $N_G[x_i]$ is played by Staller. The two cliques incident to $x_i$, say $Q_p$ and $Q_q$, together have $2t$ representatives. As $k <2t$, one of those representatives, say $w$, was played by Staller in an earlier move. Suppose that $w \in R_p$. After the move when $w$ was played, vertex $x_i \in Q_p$ remained unplayed. By the supposed strategy of Dominator, he played a vertex from $Q_p$ in his next move. As $Q_p \subseteq N_G[x_i]$, it contradicts our assumption. This contradiction proves the claim.  \hfill$(\square)$
      \medskip
      
      Claim 11.A shows that with a move $x_1, \dots , x_k$, Staller cannot claim an entire neighborhood $N_G[x_i]$, for $i \in [k]$. We observe that the same is true for each vertex $x \in N_G[x_i]$ that was played in an earlier move. Indeed, when $x$ was played, its neighbor $x_i$ remained unplayed, and Dominator could reply by choosing a neighbor of $x$. Therefore, Staller cannot claim the neighborhood $N_G[x]$ if Dominator follows the described strategy.
      
      The conclusion is that Dominator can prevent Staller from claiming a closed neighborhood of a vertex in the game and equivalently, Dominator wins the $(k,k)$-biased S-game. By Proposition~\ref{prop:W}, we conclude $\ra_k(G)\leq \ra_k'(G) \leq k$ and $\rb_k(G)\geq \rb_k'(G) \geq k+1$. 
  \end{proof}
  \begin{remark}
      M.\ Hall proved that P.\ Hall's theorem can be extended to infinite families of finite sets if the vertex degrees are finite~\cite[Theorem 5.1.2]{m-hall}. The same is true for the case of $k$-representatives. Therefore, Theorem~\ref{thm:line-graph-general} remains true if $G$ is an infinite line graph with finite vertex degrees. 
  \end{remark}
  By Theorem~\ref{thm:line-graph-general}, we can get the exact values of $\ra_1(L(H))$ and $\ra_1'(L(H))$ if $\delta(H) \ge 2.$
  \begin{corollary}
   If $G$ is the line graph of a graph $H$ with $\delta(H) \ge 2$, then $\ra_1(G)=\ra_1'(G)=1$.    
  \end{corollary}
  Theorems~\ref{thm:large} and \ref{thm:line-graph-general} together imply the following exact values for the thresholds.
  \begin{corollary} \label{cor:line-graph-large}
      Let $H$ be a graph with $\delta(H) \ge 2t$ and $G=L(H)$ its line graph with $\Delta(G)=\Delta$. Then, the following statements hold for every integer $k$ with $2 \le k \le 2t-1$.
          \begin{enumerate}[(a)]
        \item If $n(G) \ge k^{\Delta-k+1}(\Delta^2+1)$, then $\ra_k'(G) = k$ 
        and $\rb_{k-1}'(G) = k$.
         \item If $n(G) \ge (k^{\Delta-k+1}+1)(\Delta^2+1)$, then $\ra_k(G) = k$ and $\rb_{k-1}(G) = k$.
              \end{enumerate}
  \end{corollary}

  %%%%%%%%%%%%%%%%%%%%%%%%%%%%%%%%%%%%%%%%%%
  %%%%%%%%%%%%%%%%%%%%%%%%%%%%%%%%%%%%%%%%%%%%
  \section{A summary of thresholds for grids}
  \label{sec:grids}

 Concerning paths, we know that $W(P_n,1,1)=D$ holds for a path $P_n$ on $n \ge 1$ vertices; see~\cite{duchene-2020}. 
If $n \in [3]$, then $\gamma(P_n)=1$ and Proposition~\ref{prop:trivial-bounds} (iii) gives $\ra_\ell(P_n)=1$ for every positive integer $\ell$. It also implies $\rb_\ell(P_n)= \infty$ for every $\ell \ge 1$. 
For any $n \geq 4$, the two leaves of $P_n$ cannot be dominated by one vertex. Thus $W(P_n,1,2)=S$ holds for any $n \geq 4$.  We may infer that ${\rm b}_1(P_n)=2$ for any $n \geq 4.$ On the other hand $W'(P_n,1,1)=D$ holds if and only if $P_n$ has a perfect matching; see~\cite{duchene-2020}. Since $W'(P_n,1,2)=S$, it follows that ${\rm b}'_1(P_n)=2$ if $n$ is even and ${\rm b}'_1(P_n)=1$ if  $n$ is odd. Moreover, $W'(P_n,2,1)=D$, since Dominator can even achieve that no two vertices played by Staller in the $(2,1)$-biased S-game are adjacent at the end of the game. Thus, we immediately get that ${\rm a}'_1(P_n)=2$ if $n$ is odd, and ${\rm a}'_1(P_n)=1$ if $n$ is even.

It is known that Dominator wins both the $(1,1)$-biased S-game and D-game (see~\cite{dokyeesun-2024, forcan-2022+}) on every finite grid $P_m \Box P_n$ with $m \geq n \geq 2$. Consequently, 
$$\ra_1(P_m \Box P_n)=\ra_1'(P_m \Box P_n)=1,$$  
and $\rb_1 (P_m \Box P_n)\ge \rb'_1(P_m \Box P_n) \ge 2$.
Next, we consider Staller's thresholds for these grids. It is easy to observe that ${\rm b}_1(P_2 \Box P_2)=3$ and ${\rm b}'_1(P_2 \Box P_2)=2$. For $m\geq 3$ and $n \geq 2$ the result is established in the following proposition.

\begin{proposition}
If $m \ge 3$ and $n \ge 2$, then ${\rm b}_1(P_m \Box P_n)={\rm b}'_1(P_m \Box P_n)=2$.    
\end{proposition}
\begin{proof}
Let $m\geq 3, n \geq 2$ and let $V(P_m\Box P_n)=[m] \times [n]$. Consider the $(1,2)$-biased D-game played on $P_m\Box P_n$.  

\medskip 

\noindent\textit{Case 1.} Assume that $n=2$ and Dominator selects  vertex $v$ in his first move.

If $m=3$, then we may suppose, by symmetry, that $v \in \{(1,1), (2,1)\}$. If $v=(1,1)$, then Staller replies by playing $(3,1)$ and $(3,2)$. She then wins the game in her second move by playing an unplayed vertex in the middle $P_2$-layer. If $v=(2,1)$, then Staller replies by choosing $(1,2)$ and $(2,2)$. Then Dominator needs to play $(1,1)$ and Staller will win the game in her next move by claiming $(3,1)$ and $(3,2)$.

If $m \ge 4$, there is a corner such that its two neighbors are  not dominated by $v$. Without loss of generality we may assume that $(1,1),(1,2)$ and $(2,1)$ are such undominated vertices. Then Staller plays $(1,1)$ and $(1,2)$ in her first move, and she can win the game in her second move by selecting at least one of the vertices  $(2,1)$ and $(2,2)$.

\medskip 

\noindent\textit{Case 2.} Assume that $n \ge 3$ and Dominator selects vertex $v$ in his first move.

If $m \ge 4$, it remains true that after any move $v$ of Dominator, there is a corner in $P_m \Box P_n$ such that its two neighbors are  not dominated by $v$. We may assume again that $(1,1),(1,2)$ and $(2,1)$ are such undominated vertices. Staller selects vertices $(1,1)$ and $(1,2)$ in her first move. Dominator then needs to reply by playing $(2,1)$. To win the game, Staller plays $(1,3)$ and $(2,2)$ in the next move.

 Suppose now that $m=n=3$. If $v \neq (2,2)$, then Staller can apply the same strategy as for the case of $m \ge 4$. If Dominator plays  $v=(2,2)$ in his first move, then  Staller replies by choosing vertices $(1,1)$ and $(1,2)$. It forces Dominator to play $(2,1)$. With her next move, Staller can win the game by playing $(1,3)$ and $(2,3)$.
\medskip
  
We have proved that Staller can win in the $(1,2)$-biased D-game on $P_m \Box P_n$ when $m \ge 3$ and $n \ge 2$. Therefore, we have ${\rm b}_1'(P_m\Box P_n) \le {\rm b}_1(P_m\Box P_n) \le 2$. Since Staller cannot win in the $(1,1)$-biased MBD games, we may conclude ${\rm b}_1'(P_m\Box P_n) = {\rm b}_1(P_m\Box P_n)= 2$ as stated.
 \end{proof}

An \emph{infinite path} is a path $P_{\infty}=(\mathbb{Z}, E)$ where $E= \{ \{i,i+1\} \mid i \in \mathbb{Z}\}$. We can always find a perfect matching in $P_{\infty}$. This enables us to use the so-called {\em pairing strategy}, which is based on the idea that a player partitions a set of vertices into pairs, and follows the other player by selecting one vertex of each pair; see~\cite{duchene-2020}.
Applying the pairing strategy, Dominator can win both the $(1,1)$-biased D-game and S-game on $P_{\infty}$ by playing a vertex from every edge in the matching. This implies $\ra_1(P_{\infty})=\ra'_1(P_{\infty})=1$. Note that $W(P_\infty, \ell, 3)=S$ for any $\ell$, whereas $W(P_\infty, \ell, 2)=S$ only if $\ell\le 3$. 
Hence, the following thresholds for Staller hold.  

\begin{observation}
Let $\ell$ be a positive integer. Then 

\begin{align*}
    {\rm b}_{\ell}( P_{\infty})={\rm b}'_{\ell}(P_{\infty})=
                \begin{cases}
                        2; &\ell \le 3,\\
                        3; &\ell \ge 4.
                    \end{cases}
\end{align*}
\end{observation}
 
We now consider infinite grids  $P_n \Box P_{\infty}$ and $P_{\infty} \Box P_{\infty}$.  We remark that the vertex degrees are finite in $P_n \Box P_{\infty}$ and $P_{\infty} \Box P_{\infty}$.  Therefore, Staller can win the game if she can claim a closed neighborhood of a vertex, and we may say that Dominator wins the game if he has a strategy to prevent Staller from claiming a closed neighborhood. 
With this explanation, we may extend the definition of thresholds ${\rm a}_{\ell}, {\rm a}_{\ell}'$, ${\rm b}_{\ell}, {\rm b}_{\ell}'$ for infinite graphs with bounded maximum degree. For instance, if $G$ is a graph and $\ell$ a positive integer, the minimum $a$ such that Dominator can prevent Staller from claiming an entire closed neighborhood of a vertex throughout the game is denoted by ${\rm a}_{\ell}$ when Dominator starts the game on $G$.

\begin{proposition} \label{prop:infinite-grid-a1-b1}
If $n \geq 2$, the following statements hold.
\begin{enumerate}[(i)]
    \item ${\rm a}_1(P_n \Box P_{\infty})={\rm a}'_1(P_n \Box P_{\infty})={\rm a}_1(P_{\infty} \Box P_{\infty})={\rm a}'_1(P_{\infty} \Box P_{\infty})=1$,
    \item ${\rm b}_1(P_n \Box P_{\infty})={\rm b}'_1(P_n \Box P_{\infty})={\rm b}_1(P_{\infty} \Box P_{\infty})={\rm b}'_1(P_{\infty} \Box P_{\infty})=2$.  
\end{enumerate}
\end{proposition}

\begin{proof} We consider the grid  $P_n \Box P_{\infty}$ with vertex set $[n] \times \mathbb{Z}$, while $P_\infty \Box P_{\infty}$ is defined on $\mathbb{Z} \times \mathbb{Z}$.

(i) Observe that all infinite grids $P_n \Box P_{\infty}$ and $P_{\infty} \Box P_{\infty}$ have a perfect matching. Therefore,  in the (1,1)-biased MBD game, Dominator may follow the pairing strategy. For example, if Staller plays a vertex $(i,2k)$, Dominator replies with vertex $(i, 2k+1)$ and vice versa. By this strategy, Staller can never claim a closed neighborhood and hence, Dominator wins the game. 
 
(ii) 
For infinite grids, we have $\Delta(P_2 \Box P_{\infty})=3$, $\Delta(P_n \Box P_{\infty})=4$ if $n \ge 3$, and $\Delta(P_{\infty} \Box P_{\infty})=4$.  
Theorem~\ref{thm:large}\,(b) and Remark~\ref{rem:large} then imply ${\rm b}_1(P_n\Box P_{\infty}) \le 2$ and ${\rm b}_1(P_{\infty}\Box P_{\infty}) \le 2$. By part (i), we may infer the equality ${\rm b}_1(P_n\Box P_{\infty}) = {\rm b}_1(P_{\infty} \Box P_{\infty}) =2$. Since ${\rm b}'_1(G) \leq {\rm b}_1(G)$ holds for every graph $G$, we have that ${\rm b}'_1(P_n\Box P_{\infty})={\rm b}'_1(P_{\infty} \Box P_{\infty})=2$ is also true.
\end{proof}

The $(a,b)$-biased games on grids may also be interesting when $a,b \geq 2$. If $G$ is a finite grid, then Proposition~\ref{prop:trivial-bounds} implies that $\rb_a'(G) \leq 3$ for any $a \ge 1$, and $\ra_b'(G) = \infty$ for any $b \geq 3$. If $G$ is an arbitrary grid (finite or infinite), then Proposition~\ref{prop:trivial-bounds} implies that $\rb_a(G) \leq 5$ holds for any $a < \gamma(G)$. 
If $G$ is an infinite grid, then Theorem~\ref{thm:large} and Remark~\ref{rem:large} imply that $\ra_k'(G) \ge \ra_k(G) \ge k$  and $\rb_k'(G) \le \rb_k(G) \le k+1$ for every positive integer $k$. Applying also the results from Theorem~\ref{thm:k-local-dom-a1'} and Proposition~\ref{prop:infinite-grid-a1-b1}, we obtain the following estimations for an infinite grid $G$:   
$$2 \le \ra_2(G) \le \ra_2'(G) \le 4 \textrm{ and } 2 \le \rb_2'(G) \leq \rb_2(G) \leq 3;$$ 
$$3 \le \ra_3(G) \le \ra_3'(G) \le 6 \textrm{ and } 2 \le \rb_3'(G) \leq \rb_3(G) \leq 4.$$

Now, consider the finite grid $G=P_m \Box P_n$. Theorem~\ref{thm:k-local-dom-a1'} implies that $\ra_2'(G) \leq 4$. In the next result we show that Dominator has to dominate at least three vertices per move to win the $(a,2)$-biased S-game on $G$, when $m \equiv 1 \pmod 4$.

\begin{proposition}\label{p:gridsBias2}
    Let $G=P_m\Box P_n$, where $m=4k+1$ for some $k \in {\mathbb{N}}$ and $n \in \mathbb{N}$. Then $\ra_2'(G)\geq 3$.
\end{proposition}
\begin{proof}
     Let $G=P_m\Box P_n$, where $m=4k+1$ for some $k \in {\mathbb{N}}$ and let $V(G)=[m]\times [n]$. We describe Staller's strategy to win the $(2,2)$-biased S-game on $G$. In her first move Staller selects vertices $(2,1),(4,1)$. Then Dominator has to play one vertex from each of the sets $\{(1,1),(1,2)\}, \{(3,1),(3,2)\}$ to prevent Staller claiming all vertices of $N_G[(1,1)]$ or $N_G[(3,1)]$. If $k=1$, then Staller wins in her next move selecting $(5,1)$ and $(5,2)$. Otherwise, in her $i$th move, for any $i\in \{2,\ldots ,k\}$, Staller selects vertices $(4i-2,1),(4i,1)$. Note that to her $i$th move Dominator has to respond by a vertex from each of the sets $\{(4i-3,1),(4i-3,2)\}, \{(4i-1,1),(4i-1,2)\}$ to prevent Staller winning in her $(i+1)$st move. Finally, unless the game finished earlier with Staller's win, she claims all vertices from $N_G[(4k+1,1)]$ in her $k$th move, and wins. Hence $\ra_2'(G) \geq 3$.  
\end{proof}

Since $\rb_2'(G) \geq \rb_1'(G)=2$ and as Dominator does not win the $(2,2)$-biased S-game played on $P_{4k+1}\Box P_n$, we get the following.

\begin{corollary}\label{cor:gridsBias2}
    If $G=P_m\Box P_n$, where $m=4k+1$ for some $k \in {\mathbb{N}}$ and $n \in \mathbb{N}$, then $\rb_2'(G)=2$.
\end{corollary}

%%%%%%%%%%%%%%%%%%%%%%%%%%%%%%%%%%%%%%%%%%%%%%%%%%%%%%%%%%%%%%%%%%
%%%%%%%%%%%%%%%%%%%%%%%%%%%%%%%%%%%%%%%%%%%%%%%%%%%%%%%%%%%%%%%%%

\section{Star partition width and $\ra_1'(G)$} 
\label{sec:star-partition}

In this section, we introduce the concept of star partition width and prove that it serves as an upper bound for the threshold at which Dominator wins the biased games when Staller selects exactly one vertex per move. 
 
 For a positive integer $k$, we say that a graph $G$ has a {\em $k$-star partition} if there exists a partition $\{S_1,\ldots,$ $S_m\}$ of $V(G)$ such that each $S_i$ satisfies $|S_i|\ge 2$ and $G[S_i]$ contains a spanning star of order at most $k+1$.
 The {\em star partition width} of a graph $G$, $\sigma(G)$, is the smallest integer $k$ such that there exists a $k$-star partition of $G$. 
Let ${\cal S} = \{S_1,\ldots, S_\ell\}$ be a $k$-star partition of $G$. For instance, $\sigma(K_n)=1$ if $n$ is even, and $\sigma(K_n)=2$ if $n$ is odd, while $\sigma(K_{1,r})=r$. For any $i \in [\Delta(G)]$, we denote by $s_i$ the number of $i$-stars in ${\cal S}$. We say that a star partition ${\cal S}$ is \emph{lexicographically optimal} if 
 $(s_\Delta, s_{\Delta - 1}, \ldots , s_1)$ is lexicographically the smallest sequence over all star partitions of $G$. Note that the largest index $i$ such that $s_i \neq 0$ in ${\cal S}$ is $\sigma(G)$.

We remark that the so-called star partitions were considered by~\cite{andreatta-2019}, although the corresponding partitions do not restrict the size of the spanning stars as is the case with $k$-star partitions. The concept of induced star partitions as studied by~\cite{shalu} is even further away from our concept, since the stars in the corresponding partition are induced. On the other hand, there is  a close connection between $k$-star partitions and $[1,k]$-factors as studied by~\cite{heinrich-1990}. 

A {\em $[1,k]$-factor} of a graph $G$ is a spanning subgraph $F$ of $G$ such that $1 \leq \textrm{deg}_F(x) \leq k$ for every vertex $x$ in $G$. Thus, a graph $G$ has a $[1,k]$-factor if and only if it has a $k$-star partition. Hence~\cite[Corollary 3]{heinrich-1990} implies the following.
 
\begin{theorem}
\label{thm:k-star}
    Let $G$ be a graph and $k\ge 2$. Then $G$ has a $k$-star partition if and only if $i(G-X) \leq  k|X|$ for every $X \subseteq V(G)$, where $i(G-X)$ is the number of isolated vertices in $G-X$.
\end{theorem}

Note that the special case of Theorem~\ref{thm:k-star} for $k=2$ was proved by~\cite{Lovasz-1970}.
As an immediate consequence of the above theorem, we derive the formula in Corollary~\ref{cor:sigma}. A \emph{nontrivial path cover} of a graph is a spanning subgraph whose components are paths of order at least two. Clearly, a graph has a nontrivial path cover if and only if it admits a $2$-star partition.

\begin{corollary} \label{cor:sigma}
 If $G$ is a graph without a nontrivial path cover and without isolated vertices, then $$\sigma(G)=\max_{\emptyset \subset S\subset V(G)}\Big\{\Big\lceil\frac{i(G-S)}{|S|}\Big\rceil\Big\}$$    
\end{corollary}

 \begin{proposition}
 \label{prp:sigmanu}
 \begin{enumerate}[(i)]
\item For every graph $G$, the maximum number of stars in a nontrivial star partition equals $\nu(G)$. 
\item For every tree $T$, each lexicographically optimal star partition contains  $\nu(T)$ stars.
\end{enumerate}
\end{proposition}
\begin{proof}
    (i) Let $\cS$ be a nontrivial star partition of $G$. By choosing one edge from each star $S_i \in \cS$, we obtain a matching of size $|\cS|$ in $G$. Then, $|\cS| \leq \nu(G)$. Now, we consider a maximum matching $M$ in $G$ and prove that a star partition $\cS$ exists with $|\cS|=|M|=\nu(G)$. 
    By the maximality of $M$, for each edge $xy \in M$, one of the following cases holds: 
    \begin{itemize}
        \item Both $x$ and $y$ have exactly one unsaturated neighbor, and this unsaturated vertex is the common neighbor of  $x$ and $y$. We may specify either $x$ or $y$ as a center in the star partition. 
        \item Only one of $x$ and $y$ is adjacent to vertices unsaturated by $M$. If $x$ is such a vertex with an unsaturated neighbor, we specify $x$ as a center in the star partition. 
        \item Neither $x$ nor $y$ have unsaturated neighbors. Then we may choose any of them as a center. 
    \end{itemize}
     With the set of central vertices in hand, we can easily build a star partition $\cS$ of $G$ with $|\cS|=\nu(G)$. This proves that the maximum number of stars equals $\nu(G)$.

    (ii) For a tree $T$, the inequality $|\cS| \leq \nu(T)$ remains true for any star partition $\cS$. Further, if $\cS$ is lexicographically optimal, the centers can be chosen such that the leaves of the stars form an independent set in $T$. Then, the $|\cS|$ centers yield a vertex cover, and we may conclude that
    $$ \tau(T) \leq |\cS| \leq \nu(T)
    $$
    holds. As $T$ is a tree, this implies  $|\cS|= \nu(T)=\tau(T)$ for every lexicographically optimal $\cS$.   
    \end{proof}

\begin{corollary}
For a nontrivial graph $G$, we have $$\Big\lceil \frac{n(G)}{\nu(G)}\Big\rceil -1\le \sigma(G)\le \Delta(G).$$
\end{corollary}
\begin{proof}
Since any star in a graph $G$ can have at most $\Delta(G)$ leaves, the upper bound is clear. By Proposition~\ref{prp:sigmanu}(i), any nontrivial star partition of $G$ contains at most $\nu(G)$ stars. Combining this with the fact that the largest star in a lexicographically optimal star partition of $G$ has at most $\sigma(G)+1$ vertices, we infer that $(\sigma(G)+1)\nu(G)\ge n(G)$. Since $\sigma(G)$ is an integer, the result follows. 
\end{proof}
    
Next, we define the so-called {\em star digraph}, which applies to a graph $G$ with a fixed $k$-star partition ${\cal S}=\{S_1,\ldots, S_m\}$. Notably, the digraph $H_{G,\cS}=(V,A)$ has $V=\cal S$, while $S_iS_j\in A$ if a leaf of $S_i$ is adjacent in $G$ to the center of $S_j$; see Fig.~\ref{fig:starDig} for an example illustrating this construction.

\begin{figure}
\begin{center}
	\begin{tikzpicture}[]
        %edges	
	\draw (0,0) -- (1,1.5)--(2,0);
	\draw (3,0)--(4.5,1.5)--(4,0);
        \draw (5,0)--(4.5,1.5)--(6,0);
        \draw (6,0)--(7,1.5)--(7,0);
        \draw (8,0)--(7,1.5);
        \draw (1,1.5)--(4.5,1.5);

        \draw [->, line width=1pt] (13,0) -- (12.1,1.35);

	%vertices

	\draw [fill=black] (1,1.5) circle (5pt);
	\draw [fill=black] (0, 0) circle (5pt);
	\draw [fill=black] (2, 0) circle (5pt);
	\draw [fill=black] (3, 0) circle (5pt);
	\draw [fill=black] (4, 0) circle (5pt);
	\draw [fill=black] (5, 0) circle (5pt);
	\draw [fill=black] (6, 0) circle (5pt);
	\draw [fill=black] (7, 0) circle (5pt);
        \draw [fill=black] (8, 0) circle (5pt);
       \draw [fill=black] (4.5,1.5) circle (5pt);
       \draw [fill=black] (7,1.5) circle (5pt);

       \draw [fill=black] (11,0) circle (5pt);
       \draw [fill=black] (12,1.5) circle (5pt);
       \draw [fill=black] (13,0) circle (5pt);

       \draw (1,0.65) ellipse (1.3cm and 2cm);
        %\draw (1,0.65) circle (45pt);
        \draw (4,0.65) ellipse (1.5cm and 2cm);
        \draw (7,0.65) ellipse (1.4cm and 2cm);
	 
	\draw (-0.3,2) node{$S_1$}; 
        \draw (2.6,2) node{$S_2$}; 
        \draw (8.4,2) node{$S_3$}; 
        \draw (4.5,-1.8) node{A graph $G$ and a star partition ${\cal{S}}=\{S_1,S_2,S_3\}$}; 
        \draw (12,-1.8) node{$H_{G,{\cal{S}}}$}; 
        \draw (10.6,0.3) node{$S_1$}; 
        \draw (12,2) node{$S_2$}; 
        \draw (13.3,0.3) node{$S_3$};
	
	\end{tikzpicture} \\
	\caption{A graph $G$ with star partition ${\cal{S}}=\{S_1,S_2,S_3\}$ and its star digraph $H_{G,{\cal{S}}}$}  \label{fig:starDig}
	
\end{center}
\end{figure}

 \begin{lemma} \label{lem:lexico-optimal}
     Let $G$ be a nontrivial graph and let ${\cal S} = \{S_1,\ldots, S_\ell\}$ be lexicographically optimal $\sigma(G)$-star partition of $G$. Then the following statements hold for every $x, y \in V(G)$.
     \begin{enumerate}[(i)]
         \item If $x, y$ are leaves of stars in ${\cal S}$, then $xy \notin E(G)$, except possibly in the following special cases: (1) $x$ and $y$ belong to one and the same $S_i$ of which $x$ and $y$ are the only leaves, (2) $x$ and $y$ belong to distinct stars $S_i$ and $S_j$ such that $|S_i|+|S_j|\le 5$.
         \item Let $x$ be a leaf of $S_i$ and $y$ the center of $S_j$, where $i, j \in [\ell]$ and $i \neq j$. If $xy \in E(G)$, then $|S_i|\le |S_j|+1$.
         \item Let $S_{i_1}$ be a $\sigma(G)$-star in $\cS$, and let $S_{i_1}, S_{i_2}, \ldots, S_{i_r}$ be a directed path in $H_{G,\cS}$. Then $|S_{i_r}| \geq \sigma(G)$.
     \end{enumerate}
 \end{lemma}
 \begin{proof}
 (i) Let $x, y$ be the leaves of stars in ${\cal S}$, which do not satisfy the special cases, as described in the statement (i). We distinguish two cases. If $x$ and $y$ belong to the same star $S_i$ in $\cS$, then $|S_i|\ge 4$. If $xy\in E(G)$, then $\cS'$ obtained from $\cS$ by adding a new star $S'$ with $S'=\{x,y\}$ and removing $x$ and $y$ from $S_i$ yields a contradicition to $\cS$ being lexicographically optimal. Now, let $x\in S_i$ and $y\in S_j$, where $i\ne j$. Suppose $xy\in E(G)$. Note that our initial assumption gives $|S_i|+|S_j|> 5$. Let $|S_i|\ge |S_j|$. Suppose first that $|S_j|\ge 3$. Then, the star partition $\cS'$ obtained from $\cS$ by adding a new star $S'$ with $S'=\{x,y\}$ and removing $x$ from $S_i$ and $y$ from $S_j$ yields a contradicition to $\cS$ being lexicographically optimal. Finally, if $|S_j|=2$, then the star partition $\cS'$ obtained from $\cS$ by removing $x$ from $S_i$ and adding it to $S_j$ yields a contradicition to $\cS$ being lexicographically optimal.

 (ii) Suppose a leaf $x$ of $S_i$ and the center $y$ of $S_j$, where $i \neq j$, are adjacent, and suppose to the contrary that $|S_i|\ge |S_j|+2$. Then the star partition $\cS'$ obtained from $\cS$ by removing $x$ from $S_i$ and adding it to $S_j$ yields a contradiction to $\cS$ being lexicographically optimal.

 (iii) Let $S_{i_1}$ be a $\sigma(G)$-star in $\cS$, and let $S_{i_1}, S_{i_2}, \ldots, S_{i_r}$ be a directed path in $H_{G,\cS}$. For every $j\in [r]$, denote by $x_j$ the center of the star $S_{i_j}$, and for every $j\in [r-1]$ denote by $y_j$ a leaf of $S_j$, which is adjacent to $x_{j+1}$. Suppose to the contrary, $|S_{i_r}|<\sigma(G)$. Then, the star partition $\cS'$ obtained from $\cS$ by removing $y_j$ from $S_{i_j}$ and adding it to $S_{i_{j+1}}$ for all $j\in [r-1]$ yields a contradiction to $\cS$ being lexicographically optimal.
  \end{proof}

We use the above lemma to prove that the star partition width of $G$ is an upper bound for ${\rm a}_1'(G)$ and get the exact value of ${\rm a}_1'(T)$ for any tree $T$.

   \begin{theorem} \label{thm:sigma-a1'}
     If $G$ is a nontrivial graph, then ${\rm a}_1'(G) \leq \sigma(G)$. In particular, ${\rm a}_1'(T) = \sigma(T)$ holds for every nontrivial tree $T$.
 \end{theorem}
 \begin{proof}
 Let $G$ be a nontrivial graph and $k= \sigma(G)$. Then there exists a $k$-star partition of $G$, say $\{S_1,\ldots, S_\ell\}$, such that each $S_i$ contains a spanning star with at most $k$ leaves. Assume that Staller plays a vertex in a star $S_i$ for $i \in [\ell]$. Then the strategy for Dominator is to select all the remaining vertices in $S_i$ and arbitrary $k-|S_i|+1$ vertices. By this strategy, Staller cannot claim a closed neighborhood of a vertex in $G$, that is $W'(G,k,1)=D$. We may then conclude ${\rm a}_1'(G) \leq \sigma(G)$.

 Assume next that $T$ is a tree. If $\sigma(T)=1$, then $\ra'_1(T) \le \sigma(T)$ implies that equality holds.  
 
 If $\sigma(T)\ge 2$, we set $a=\sigma(T)-1$ and prove that $W'(T,a,1)=S$. We proceed by induction on the matching number of $T$. If $\nu(T)=1$, then $T$ is a star, and the statement clearly holds (Staller plays the center, Dominator selects $a$ leaves, and Staller wins by selecting an unplayed leaf). Now, let $T$ be a tree with $\nu(T)>1$.     
    Let ${\cal S} = \{S_1,\ldots, S_\nu\}$ be a lexicographically optimal $\sigma(G)$-star partition of $G$. Let $S_{i_0}\in\cS$ be a $\sigma(G)$-star in $\cal S$. Let $S_{i_0},\ldots, S_{i_t}$ be a maximal directed path in the star digraph $H_{T,\cS}$. Therefore, since $T$ is a tree, $S_{i_t}$ has no out-edges in $H_{T,\cS}$. By Lemma~\ref{lem:lexico-optimal}(iii), we have $|S_{i_j}|\ge \sigma(T)$. 

    The strategy of Staller to win the game in $T$ is as follows. In her first move, Staller selects the center of $S_{i_t}$, after which Dominator needs to select all the leaves of $S_{i_t}$ to avoid an immediate defeat in the next move of Staller. (In particular, 
    if $|S_{i_j}|=\sigma(T)+1$, the game is over.) Now, consider the tree $T'$, which is the connected component of $T-S_{i_t}$, which contains $S_{i_0}$. Note that $\nu(T')<\nu(T)$, and that $\sigma(T')=\sigma(T)$ (the latter holds, because $\cS$ is a lexicographically optimal star partition of $T$). Therefore, by the induction hypothesis, Staller has a winning strategy in $T'$. In the rest of the game, Staller plays only in $T'$ by using her winning strategy.   

    The inequality ${\rm a}_1'(T) \leq \sigma(T)$ and $W'(T,\sigma(T)-1,1)=S$ together imply ${\rm a}_1'(T) = \sigma(T)$ for every tree with $\sigma(T) \ge 2$.
 \end{proof}

We show that there exist graphs $G$
 with ${\rm a}_1'(G) <\sigma(G)$, and the difference can be arbitrarily large. Consider a complete bipartite graph $K_{2,2m}$, with $u,v$ being the vertices of degree $2m$. Dominator can win the $(1,1)$-biased S-game by playing $u$ or $v$ as his first move and an arbitrary vertex in his second move. Hence ${\rm a}'_1(K_{2,2m})=1$. As $\sigma(K_{2,2m})=m$, we deduce that for any $n \in \mathbb{N}$ there exists a graph $G$ with ${\rm a}_1'(G) = \sigma(G)-n$. 
 
 \medskip

 We remark that manuscript of~\cite{bagdas-2024} characterizes trees with $W'(T,k,1)=D$. It can be proved directly that the definition of a $k$-good tree in~\cite{bagdas-2024} corresponds to the condition $\sigma(T) \le k$ if $T$ is a tree and therefore, \cite[Theorem~5]{bagdas-2024} directly implies our result on trees. Anyway, we provided a different proof based on the star partition width, a concept that could be of independent interest.
 \medskip

 In Theorems~\ref{thm:k-local-dom-a1'} and \ref{thm:sigma-a1'} we established two upper bounds for $\ra_1'(G)$ that are incomparable. For a tree $G$, $\widetilde{\gamma}(G) = \Delta(G)$ while $\sigma(G)$ may be much smaller than $\Delta(G)$. On the other hand, if $G$ is the complete bipartite graph $K_{2, 2m}$, then $\widetilde{\gamma}(G) =1$ and $\sigma(G)=m$. 

\section{Concluding remarks}

Note that the tightness of the bound in Proposition~\ref{prop:trivial-bounds}(iv), namely, ${\rm a}'_b(G) \leq b \cdot\Delta(G)$ when $b \le \delta(G)$, was provided only in the case $b=1$, where already odd paths attain the bound. We do not know if for larger values of $b$ the bound is also tight and pose it as an open problem.

\begin{problem}
Is it true that the bound ${\rm a}'_b(G) \leq b \cdot\Delta(G)$ is tight also when $1<b\le \delta$. 
\end{problem}

From Proposition~\ref{p:gridsBias2} and the earlier remark, we derive that $3\le  \ra_2'(P_{4k+1}\Box P_n)\leq 4$, and it would be interesting to determine which of the two bounds is the exact value. For grids $P_m\Box P_n$, where none of the integers $n$ and $m$ is congruent to $1$ modulo $4$, the currently known bounds are further apart:  $1\le  \ra_2'(P_m\Box P_n)\leq 4$. (Remark that, by Theorem~\ref{thm:large}, the lower bound can be replaced with $2$ when the order of $P_m\Box P_n$ is at least $136$.)  Based on the discussion, we pose the following problem.

\begin{problem}
\label{prob:grids}
Determine the values $\ra_2'(P_m\Box P_n)$ for all $m,n\in \mathbb{N}$, and the values $\rb_2'(P_m\Box P_n)$ when none of the integers $m$ and $n$ is congruent to $1$ modulo $4$.
\end{problem}

Problem~\ref{prob:grids} can also be extended to infinite grids, which might be even more difficult. In addition, determining the values of $\ra_2(P_m\Box P_n)$ and $\rb_2(P_m\Box P_n)$ for all $m,n\in \mathbb{N}$ is also a challenging problem.   

In Section~\ref{sec:upper-bound-a1}, we obtained as a corollary the bound ${\rm a}_1'(G) \leq k-1$, which holds for any induced $K_{1,k}$-free graph $G$ with $\delta(G) \ge 1$, and proved its sharpness in the case when $k=3$ (that is, in claw-free graphs). We do not know if the sharpness holds for larger $k$.

\begin{problem}
\label{prob:a1}
Are there any $K_{1,k}$-free graphs $G$, where $k>3$ such that ${\rm a}_1'(G)=k-1$?
\end{problem}

\acknowledgements
\label{sec:ack}
We wish to thank the anonymous reviewers for several useful remarks. %This work was supported in part by the Slovenian Research and Innovation Agency (ARIS) under the grants P1-0297, N1-0285, N1-0355, and J1-4008. 
The authors also thank the Institute of Mathematics, Physics and Mechanics for hosting the
2nd Workshop on Games on Graphs in June 2024 during which this work was initiated.

\nocite{*}
\bibliographystyle{abbrvnat}
% use the following instead if you encounter problems 
%\bibliographystyle{alpha}
\bibliography{dmtcs-MBD-final}
\label{sec:biblio}

\end{document}